\documentstyle[12pt]{article}
\newcommand{\beq}{\begin{equation}}
\newcommand{\eeq}{\end{equation}}
\newcommand{\barr}{\begin{array}}
\newcommand{\earr}{\end{array}}
\newtheorem{theorem}{Theorem}
\newtheorem{lemma}{Lemma}
\newtheorem{corollary}{Corollary}
\newtheorem{remark}{Remark}
\newtheorem{definition}{Definition}
\newtheorem{example}{Example}
\newtheorem{statement}{Statement}
\newtheorem{conjecture}{Conjecture}
\newcounter{three}
\input amssym.def
\begin{document}
\begin{center}
{\Huge Conformal invariant functionals of immersions of tori into ${\bf R}^3$.}

\vspace{1cm}

{\bf P.G.Grinevich}\footnote{This work was fulfilled during the
author's visit to the Freie Universit\"at, Berlin, Germany, which 
was supported by the Humboldt-Foundation. This work was also partially
supported by Russian Foundation for Basic Research, grant No 95-01-755.}

\medskip

{\bf M.U.Schmidt}\footnote{Supported by the DFG, SFB 288
``Differentialgeometrie und  Quiantenphysik''.}

\end{center}

\vspace{0.5cm}

${\vphantom k}^{\rm1}$ L.D.Landau Institute for Theoretical Physics, 
Kosygina 2, Moscow, 177940, Russia.\\ e-mail: pgg@landau.ac.ru

${\vphantom k}^{\rm2}$ Freie Universit\"at, Berlin, Arnimallee 14, D-14195,
Berlin, Germany.\\e-mail: mschmidt@omega.physik.fu-berlin.de

\begin{abstract}
We show, that higher analogs of the Willmore functional, defined on
the space of immersions $M^2\rightarrow{\Bbb R}^3$, where $M^2$ is a
two-dimensional torus, ${\Bbb R}^3$ is the 3-dimensional Euclidean
space are invariant under conformal transformations of ${\Bbb R}^3$. 
This hypothesis was formulated recently by I.~A.~Taimanov. 

Higher analogs of the Willmore functional are defined in terms of the
Modified Novikov-Veselov hierarchy. This soliton hierarchy is
associated with the zero-energy scattering problem for the
two-dimensional Dirac operator.
\end{abstract}

\section{Introduction}
\label{sec:intr}
To start with, we would like to recall the following interesting 
fact from the theory of 2-dimensional surfaces in ${\Bbb R}^3$ 
(see \cite{W93}, p. 110 and references therein). 
Let $X:M^2\rightarrow{\Bbb R}^3$ be a smooth immersion of a compact 
orientable surface $M^2$ into the Euclidean space ${\Bbb R}^3$ 
(i.e. a smooth map from $M^2$ to ${\Bbb R}^3$ such that its Jacobi
matrix is non-degenerate everywhere on $M^2$, but the image is allowed
to self-intersect). Let $T$ be the Willmore functional
\beq
T=\int_{M^2} H^2 dS,
\label{willmore}
\eeq
where $H$ denotes the mean curvature, $dS$ is the volume element 
on $M^2$ generated by the immersion. Then $T$ is invariant under
conformal transformations of ${\Bbb R^3}$.

It is natural to pose the problem of constructing other conformal 
invariant functional of immersions.

It is well-known, that many constructions from the
soliton theory have natural analogs in geometry and vise versa. 
In particular, important information how to study immersions of 
2-dimensional surfaces into ${\Bbb R}^3$ using soliton methods 
can be found in the article \cite{Bob94} by A.~I.~Bobenko. Any immersed 
surface possesses (at least locally) a conformal coordinate system 
(see for example \cite{DFN84} p. 110), i.e. a coordinate system
such that $ds^2=f(z,\bar z)dzd\bar z$. In conformal coordinates
this immersion can be locally represented by the Generalized
Weierstrass Formulas (see Section \ref{sec:weier} below) 
and the potential $U(z,\bar z)$ 
is uniquely defined. In \cite{Taim95} it was shown,
that any analytic immersion of a compact orientable 2-dimensional
manifold into ${\Bbb R}^3$ can be globally represented by these
formulas.

The Generalized Weierstrass Formulas are based on the zero-energy 
eigenfunctions of the two-dimensional Dirac operator with a real
potential $U(z,\bar z)$. The zero-energy spectral problem for this
operator arose in the soliton theory as an auxiliary linear problem
for the hierarchy of Modified Novikov-Veselov Equations (MNV) (see
\cite{Bog87}). These nonlinear integrable equations with 2 spatial
variables, introduced by L.~V.~Bogdanov ,have infinitely many
conservation laws.  

In \cite{KT93} B.~G.~Konopelchenko and I.~A.~Taimanov  showed
that the quadratic MNV conservation law 
\beq 
H_1=4\int\int U^2(z,\bar z) dx\wedge dy,\ z=x+iy
\label{energy}
\eeq
coincides with the Willmore functional $T$. In \cite{Taim96} 
I.~A.~Taimanov formulated a hypothesis that all higher MNV
conservations laws also generate functionals on immersions of closed
orientable 2-dimensional surfaces into ${\Bbb R}^3$, invariant under
conformal transformations of ${\Bbb R}^3$. He also did some numerical 
experiments with surfaces of revolution, confirming this assumption. 

During his visit to the Freie Universit\"at, Berlin in
September--October 1996, I.~A.~Taimanov
attracted the authors attention to this problem. In the present text
we prove this hypothesis for immersions of tori into ${\Bbb R}^3$. 

Of course, it would be natural to extend this result to immersions of
higher genus surfaces into ${\Bbb R}^3$. But here we meet the following
problem. In contrast with $H_1$ the higher conservation laws are
non-local in terms of the potential. For double-periodic potentials
they are defined in terms of the zero-energy dispersion curve 
(the Riemann surface of the zero-energy Bloch function). If the
genus is greater than 1, we have to study modular invariant Dirac
operators in the Lobachevskian plane with a non-abelian group of
translations. The corresponding Bloch theory has not been constructed
until now, thus we could not define the higher conservation laws. 
It would be interesting to develop the corresponding Bloch theory, 
but this problem looks rather non-trivial. Thus we restrict ourself 
to the genus 1 case only, where all integrals of motion are well-defined. 

The idea of our proof is the following. We show, that infinitesimal 
conformal transformations of ${\Bbb R}^3$ correspond to infinitesimal 
Darboux transformations of the Dirac operator (infinitesimal dressings 
with degenerate kernels). It is convenient to express these deformations
in terms of Cauchy-Baker-Akhiezer kernels, introduced by A.~Yu.~Orlov
and one of the authors in \cite{GO89} (see also the review \cite{GOS93}).
From the explicit formulas for such deformations it follows, that the
deformed zero-energy Bloch function is meromorphic on the same Riemann
surface as the original one. Recalling, that this Riemann surface
completely determines all conservation laws, we complete the proof.

\begin{remark} An alternative proof of the theorem, 
  namely that conformal transformations do not change the
  zero-energy Bloch variety, was obtained by U.~Pinkall (private
  communication). He calculated the
  action of finite conformal transformations on the Bloch function
  using a quaternionic representation of the Generalized Weierstrass
  formulas as suggested by G.~Kamberov, F.~Pedit and U.~Pinkall in \cite{KPP}.
\end{remark}
\begin{remark} Finite Darboux
transformations for the Dirac operator are discussed in the book
by V.~B.~Matveev and M.~A.~Salle \cite{MS91}, Laplace transformations
for the Dirac operator are discussed in the paper by E.~V.~Ferapontov 
\cite{Fer96}.
\end{remark}
The authors are grateful to I.~A.~Taimanov for interesting discussions.

\section{Generalized Weierstrass construction}
\label{sec:weier}
Let $L$ be the two-dimensional Dirac operator
\beq
L=\left[\barr{cc} \partial_z & -U(z,\bar z) \\ U(z,\bar z) & 
\partial_{\bar z}
\earr\right]              
\label{diracop}
\eeq
with a real potential $U(z,\bar z)$. Let $\vec\Psi(z,\bar z)$ be a
zero-energy solution of the Dirac equation
\beq
L\vec\Psi(z,\bar z)=0,\
\vec\Psi(z,\bar z)=\left(\barr{c}\Psi_1(z,\bar z)\\
\Psi_2(z,\bar z)\earr\right).
\label{diraceq}
\eeq
Then the {\bf Generalized Weierstrass Formulas} (see \cite{Taim95} and
refs. therein)
\beq
\barr{c}
X_1(z,\bar z)+iX_2(z,\bar z)=C_1+iC_2+
i\int\limits_{z_0}^{z} \left(\bar\Psi_1^2(z',\bar z') dz'- 
\bar\Psi_2^2(z',\bar z') d\bar z' \right)
\\
\\
X_1(z,\bar z)-iX_2(z,\bar z)=C_1-iC_2+
i\int\limits_{z_0}^{z} \left(\Psi_2^2(z',\bar z') dz'- 
\Psi_1^2(z',\bar z') d\bar z' \right)
\\
\\
X_3(z,\bar z)=C_3
-\int\limits_{z_0}^{z} \left(\Psi_2(z',\bar z')\bar\Psi_1(z',\bar z') dz'+ 
\Psi_1(z',\bar z')\bar\Psi_2(z',\bar z') d\bar z' \right)
\earr
\label{weier}
\eeq
defines a map of the plane ${\Bbb R}^2$ to the Euclidean space 
${\Bbb R}^3$. In (\ref{weier}) $z_0$ is a fixed point in the $z$-plane
and the integrals are taken over some path connecting the points $z_0$ and 
$z$. From (\ref{diraceq}) it follows, that the integrands in (\ref{weier})
are closed forms, thus the map does not depend on a specific choice of
the path. Here $C_1$, $C_2$, $C_3$ are arbitrary real constants.

The formulas (\ref{weier}) are equivalent to:
\beq
\barr{c}
d\left[\sigma^2 X_1+\sigma^1 X_2-\sigma^3 X_3\right] =
\\ \\
=\left[\barr{cc} \bar\Psi_1(z,\bar z)&\bar\Psi_2(z,\bar z) \\ 
-\Psi_2(z,\bar z) & \Psi_1(z,\bar z) \earr\right]
\left[\barr{cc} 0&dz\\ d\bar z&0 \earr\right]
\left[\barr{cc} \Psi_1(z,\bar z) & -\bar\Psi_2(z,\bar z) \\ 
\Psi_2(z,\bar z) & \bar\Psi_1(z,\bar z) \earr\right]
\earr
\label{weier2}
\eeq
where
$\sigma^1$, $\sigma^2$, $\sigma^3$ are the standard Dirac matrices
\beq
\sigma^1=\left[\barr{cc} 0&1\\ 1&0 \earr\right],\  
\sigma^2=\left[\barr{cc} 0&-i\\ i&0 \earr\right],\
\sigma^3=\left[\barr{cc} 1&0\\ 0&-1 \earr\right].
\label{diracmatrices}
\eeq

The Generalized Weierstrass Map is conformal, i.e. the metric $d\tilde s^2$
on ${\Bbb R}^2$ induced by this map is proportional to 
the standard one: $d\tilde s^2=g(z,\bar z)dzd\bar z$. 

Assume that we have a map of a 2-dimensional torus into ${\Bbb R}^3$. 
Then the corresponding potentials $U(z,\bar z)$ is periodic
\beq
U(z+\bar T_1,\bar z+T_1)=U(z+T_2,\bar z+\bar T_2)=
U(z,\bar z). 
\label{u}
\eeq
Also the eigenfunction $\vec\Psi(z,\bar z)$ is
periodic or anti-periodic, i.e.
\beq
\vec\Psi(z+T_1,\bar z+\bar T_1)={\cal W}_1\vec\Psi(z,\bar z),\ 
\Psi(z+T_2,\bar z+\bar T_2)={\cal W}_2\vec\Psi(z,\bar z), \ \ 
{\cal W}_1^2={\cal W}_2^2=1.
\label{persol}
\eeq

The coordinate $z$ is defined up to linear transformations
$z\rightarrow az+b$, $a,b\in{\Bbb C}$, $a\ne 0$. Thus without loss 
of generality we may assume 
\beq
T_1=1,\ T_2=\tau,\ \hbox{Im}\,\tau>0.
\eeq

Conditions (\ref{u}) and (\ref{persol}) are necessary, but,  
of course, not sufficient for periodicity of the Generalized 
Weierstrass map. Necessary and sufficient conditions for periodicity
can be formulated in terms of the Bloch variety. They are obtained in
a forthcoming paper by I.~A.~Taimanov and one of the authors
(M.S.). We do not use these conditions in our text, thus we will not
discuss them in further details

\section{Bloch function and Bloch variety.}
\label{sec:bloch}

In this Section we assume that $U(z,\bar z)$ is real, smooth, and
double-periodic (\ref{u}). With any such potential we associate a
one-dimensional subvariety $\Gamma$ in the two-dimensional complex
space $({\Bbb C}\backslash 0)^2$. 

The first object we need is the {\bf Bloch function}. By definition,
the Bloch functions $\vec\psi(w_1,w_2,z,\bar z)$ are
quasiperiodic solutions of the Dirac equation (\ref{diraceq}) 
with the following periodicity properties:
\beq
\barr{c}
\vec\psi(w_1,w_2,z+1,\bar z+1)=w_1\vec\psi(w_1,w_2,z,\bar z),
\\ 
\vec\psi(w_1,w_2,z+\tau,\bar z+\bar\tau)=w_2\vec\psi(w_1,w_2,z,\bar z).
\earr
\label{blochsol}
\eeq

The pairs of multipliers $w_1$, $w_2$ possessing at least one non-zero
Bloch solution form a complex one-dimensional subvariety
$\Gamma\in({\Bbb C}\backslash0)\times({\Bbb C}\backslash 0)$ 
(see \cite{Kuch93}). This variety is called the {\bf Bloch 
variety} or the {\bf zero-energy dispersion curve}. For a generic potential 
$U(z,\bar z)$ the genus of $\Gamma$ is infinite.

The Bloch functions form a one-dimensional holomorphic bundle over
$\Gamma$ (it is shown below, that for generic $\lambda\in\Gamma$ a
Bloch solution is unique up to normalization). 
It is convenient to fix a section of this bundle $\vec\psi(\lambda,z,\bar z)$,
by assuming
\beq
\left.\psi_1(\lambda,z,\bar z)+\psi_2(\lambda,z,\bar z)\right|_{z=z_1}=1,
\label{norm}
\eeq
where $z_1$ is an arbitrary fixed point. 

The logarithms of the multipliers $w_1(\lambda)$, $w_2(\lambda)$ 
\beq
p_1(\lambda)=\frac{1}{i} \ln w_1(\lambda),\ \ 
p_2(\lambda)=\frac{1}{i \left|\tau\right|} \ln w_2(\lambda),
\label{p's}
\eeq
are called {\bf quasimomentum functions}. Of course, they have non-trivial
increments while going along cycles in $\Gamma$, and they are defined up to
adding $2\pi n_1$, $2\pi n_2/\left|\tau\right|$ respectively, where 
$n_1$ and $n_2$ are some integers. Thus the functions 
$\hbox{Im}\,p_1(\lambda)$, $\hbox{Im}\,p_2(\lambda)$ are single-valued
in $\Gamma$. The {\bf differentials of the quasimomentum functions} 
\beq
dp_1(\lambda)= \frac{\partial}{\partial \lambda} p_1(\lambda)d\lambda, \ \ 
dp_2(\lambda)= \frac{\partial}{\partial \lambda} p_2(\lambda)d\lambda,
\label{dp's}
\eeq
are single-valued and holomorphic on the finite part of $\Gamma$.

In our text the Dirac operator (\ref{diracop}) is symmetric and the
potential $U(z,\bar{z})$ is real. Let us show, that the corresponding
Bloch variety $\Gamma$ possesses ${\Bbb Z}_2\times{\Bbb Z}_2$ as a group of
symmetries. An analogous statement for the fixed-energy Bloch variety
corresponding to a two-dimensional self-adjoint Schr\"odinger operator 
was proved in \cite {Kri89}. The proof from \cite{Kri89} may be applied
to (\ref{diracop}) after a minimal modification.

The operator (\ref{diracop}) with real potential $U(z,\bar z)$ has
the following symmetry. If $\vec\psi(w_1,w_2,z,\bar z)$ is a Bloch
solution of (\ref{diraceq}), then the function
\beq
\vec\psi^\dag(w_1,w_2,z,\bar z)=
\left ( \barr{c} \bar\psi_2(w_1,w_2,z,\bar z) \\
        -\bar\psi_1(w_1,w_2,z,\bar z) \earr \right )
\label{psidag}
\eeq
is also a Bloch solution of (\ref{diraceq}) with multipliers 
$\bar w_1$ and $\bar w_2$ respectively. Thus the surface $\Gamma$
possesses an antiholomorphic involution (we denote it by $\sigma\theta$
for historical reasons)
\beq
\sigma\theta: \Gamma\rightarrow \Gamma, \ 
\sigma\theta: (w_1,w_2) \rightarrow (\bar w_1, \bar w_2), 
\label{sigmatheta}
\eeq
and
\beq
\vec\psi^\dag(\lambda,z,\bar z)=
n^{-1}(\lambda)\vec\psi(\sigma\theta(\lambda),z,\bar z),
\label{n}
\eeq
where $n(\lambda)$ is a scalar function, meromorphic in $\lambda$ 
and independent on $z$, $\bar z$.

It is less trivial to see that the surface $\Gamma$ possesses a
holomorphic involution:
\beq
\sigma: \Gamma\rightarrow\Gamma,\ 
\sigma: (w_1,w_2) \rightarrow (w_1^{-1}, w_2^{-1}).
\label{sigma}
\eeq

To prove it, let us fix a generic point $\lambda\in\Gamma$. Denote by 
${\cal L}_{w_1 w_2}$ the Banach space of all locally square-integrable
two-component complex-valued vector-functions on ${\Bbb R}^2$ 
with the periodicity properties
(\ref{blochsol}). The space ${\cal L}_{w_1^{-1} w_2^{-1}}$ is naturally
dual to ${\cal L}_{w_1 w_2}$. Namely, if 
$\vec f(z,\bar z)\in{\cal L}_{w_1 w_2}$ and 
$\vec g(z,\bar z)\in{\cal L}_{w_1^{-1} w_2^{-1}}$, then we define a
scalar product by
\beq
\barr{r}
<f,g>=\int\limits_0^1\int\limits_0^1 dt_1\wedge dt_2 \left[
f_1(t_1+\tau t_2,t_1+\bar\tau t_2)g_1(t_1+\tau t_2,t_1+\bar\tau t_2)+\right.
\\ 
\left. +f_2(t_1+\tau t_2,t_1+\bar\tau t_2)g_2(t_1+\tau t_2,t_1+\bar\tau t_2)
\right] . \earr
\label{dual}
\eeq
Let $\vec f^{(0)}=\vec\psi(w_1,w_2,z,\bar z)$, $\vec f^{(1)}$, \ldots,
$\vec f^{(n)}$, \ldots be the Jordan basis for the Dirac operator $L$
in the space ${\cal L}_{w_1 w_2}$. Also let $\vec g^{(0)}$, $\vec
g^{(1)}$, \ldots, 
$\vec g^{(n)}$, \ldots be the dual basis in ${\cal L}_{w_1^{-1}w_2^{-1}}$. 
The functions $\vec g^{(n)}$ form a Jordan basis for the transposed 
operator $L^T$. But $L$ is symmetric with respect to this scalar product
thus it has a zero eigenfunction in the space 
${\cal L}_{w_1^{-1}w_2^{-1}}$. Hence if 
$(w_1,w_2)\in\Gamma$, then $\sigma (w_1,w_2)=(w_1^{-1}, w_2^{-1})\in\Gamma$.

One of the main properties of the Bloch variety $\Gamma$ is the
following: $\Gamma$ may be treated as a {\bf complete set of integrals
of motion} for the {\bf Modified Novikov-Veselov hierarchy}.

Indeed, consider the space of all real-valued smooth double-periodic functions
on ${\Bbb C}^1={\Bbb R}^2$ with a fixed pair of periods 1 and $\tau$. 
The Modified Novikov-Veselov hierarchy (MNV) (see
Section~\ref{sec:cons} below) defines an infinite collection of flows
on this space 
\beq
\frac{\partial U(z,\bar{z},t_{2n+1})}{\partial t_{2n+1}}= K_{2n+1}[U]+
\bar K_{2n+1}[U],
\eeq
\beq
\frac{\partial U(z,\bar{z},\tilde t_{2n+1})}{\partial t_{2n+1}}= 
i\left(\bar K_{2n+1}[U]- K_{2n+1}[U]\right),
\eeq
where $ K_{2n+1}[U]$ is some integro-differential operator in $z$, $\bar{z}$.
Here $t_{2n+1}$, $\tilde t_{2n+1}$ are parameters of these flows.

\begin{statement}
\label{st1}
Let $U(z,\bar{z},t_{2n+1})$ be a solution of one of the MNV
equations; let $L(t_{2n+1})$ be the corresponding two-dimensional Dirac 
operator (\ref{diracop}) depending on an extra parameter $t_{2n+1}$; 
let $\Gamma(t_{2n+1})$ be the corresponding family of Bloch varieties.

Then $\Gamma(t_{2n+1})=\Gamma$ does not depend on the MNV time $t_{2n+1}$.
\end{statement}

\begin{remark}
  Here and below we use the following notational convention.
  If we have a complete proof of a mathematical result we
  call it {\bf Theorem} or {\bf Lemma}. If we do not have a complete
  strict proof yet we use the word {\bf Statement}. 
\end{remark}

An analogous statement is well-known for soliton systems with
one spatial variable. In Section~\ref{sec:cons} we prove this fact at
least for algebraic-geometrical potentials, corresponding to varieties 
$\Gamma$ of finite genus. (Sometimes such potentials are called
finite-gap potentials). It is rather clear, that
our proof can be extended to all smooth potentials, but to do such
extension strictly we need more detailed information about analytic
behavior of the Bloch functions near infinity in the momentum space,
than we have now. An appropriate analytic lemma for the one-dimensional
Dirac operator, corresponding to the surfaces of revolution, was
proved by one of the authors (M.S.) in \cite{Sch95}.  

There exists a different way (may be a more natural one) to get a strict 
proof of Statement~\ref{st1}. It would be interesting to prove the
following approximation property:

\begin{conjecture}
Any smooth potential can be approximated by the algebraic-geometrical
ones with the same periods.
\end{conjecture}

From such a result it would follow, that we can restrict ourself to
the algebraic-geometrical potentials in our calculations.

We have a map $U(z,\bar{z})\rightarrow\Gamma[U]$ from the space of 
double-periodic real smooth
functions to the space of complex subvarieties in 
$({\Bbb C}\backslash 0)^2$, which is 
invariant under the whole MNV hierarchy. This map generates an
infinite family of MNV ``normal'' conservation laws. Namely, let $w_1$
be a generic point in ${\Bbb C}\backslash0$. Than we have an infinite
collection of numbers $w_2^{(k)}[U]$ such that
$(w_1,w_2^{(k)}[U])\in\Gamma[U]$. From Statement~\ref{st1} it follows
that these functionals $w_2^{(k)}[U]$ are {\bf conservation laws} of  
the whole {\bf Modified Novikov-Veselov} hierarchy.

The functionals $w_2^{(k)}[U]$ are essentially nonlocal. In 
Section~\ref{sec:cons} we show that under the same assumptions as in
Statement~\ref{st1} we can expand these functionals in some asymptotic
series near infinity and the expansion coefficients give us 
the standard ``quasi-local'' conservation laws. 

\section{Conformal transformations of the Euclidean space ${\bf R}^3$, 
and MNV integrals of motion.}
\label{sec:conf}

In the previous Section we have associated with any double-periodic
smooth real potential a Bloch variety $\Gamma[U]$. The map is constant
on the trajectories of the MNV hierarchy. In this Section we associate 
with any immersion of a torus into ${\Bbb R}^3$ a Bloch variety
$\Gamma$ and show, that $\Gamma$ is invariant under conformal 
transformations of ${\Bbb R}^3$.

Let $M^2$ be a torus with a fixed basis of cycles $a$, $b$. Let
$X:M^2\rightarrow{\Bbb R}^3$ be a smooth immersion of $M^2$
into the Euclidean space. The standard metric on ${\Bbb R}^3$ induces
a conformal structure on $M^2$. Let $z$ be a conformal global
coordinate on the universal covering space of $M^2$. The coordinate $z$ is
defined uniquely up to affine transformations $z\rightarrow cz+d$.
If we assume, that the shift of $M^2$ along the $a$-cycle corresponds to 
the shift $z\rightarrow z+1$ then the coordinate $z$ is defined
uniquely up to shifts
\beq
z\rightarrow z+d.
\label{shifts_z}
\eeq

The immersion $X$ and coordinate $z$ define a potential $U(z,\bar z)$
therefore also a Bloch variety $\Gamma[U]$. $\Gamma[U]$ is invariant
under the shifts (\ref{shifts_z}) thus it is completely
determined by the immersion $X$ and the cycles $a,b$, and we may write 
$\Gamma[X,a,b]$. 

Conformal transformations of the Euclidean space ${\Bbb R}^3$ do not
affect the conformal structure of $M^2$, thus they leave the
coordinate $z$ invariant up to the shifts (\ref{shifts_z}). Without
loss of generality we shall assume that conformal transformations of 
${\Bbb R}^3$ do not change $z$.

Now we are in position to formulate and prove our main result:

\begin{theorem} 
\label{main_theorem}  
  Let $X:M^2\rightarrow{\Bbb R}^3$ be an immersion of a
  torus with a fixed basis of cycles $a$, $b$ into the Euclidean
  space; let $\Gamma[X,a,b]$, be the corresponding Bloch variety.
  Then $\Gamma[X,a,b]$
  is invariant under conformal transformations of ${\Bbb R}^3$.
\end{theorem}

Proof of the Theorem:

{\bf Step 1}: To start with let us recall the well-known facts about the group
of conformal transformations of the standard Euclidean metric on 
${\Bbb R}^3$ (or on the sphere $S^3$) (see for example \cite{DFN84}).
This group is generated by the following transformations:

\begin{enumerate}
\item Translations $X_i\rightarrow X_i+(X_0)_i$.
\item Rotations $\vec X \rightarrow A \vec X$, $A\in SO(3)$.
\item Dilations $\vec X \rightarrow k \vec X$, $k\in {\Bbb R}$
\item Inversions 
\beq
X_i\rightarrow \frac{X_i-(X_0)_i}{<\vec X-\vec X_0,\vec X-\vec X_0>}.
\label{inversions}
\eeq
\item Reflections 
\beq
\vec X \rightarrow \vec X - 2\vec v <\vec v,\vec X>,\ <\vec v,\vec v>=1.
\eeq
\end{enumerate}

The connected component of the identity of this group is isomorphic 
to $SO(1,4)$ (see \cite{DFN84}, page 143). The corresponding Lie algebra
is generated by the following basis of infinitesimal transformations:
\begin{enumerate}
\item Translations $P_a$: $\delta X_i=\delta_{ia}$.
\item Rotations $\Omega_{ab},\ a<b$: 
$\delta X_i= \delta_{ib}X_a-\delta_{ia}X_b $.
\item Dilation $D$: $\delta X_i = X_i$
\item Inversions $K_a$
\beq
\delta X_i= 2 X_i X_a - \delta_{ia} \sum_{j=1}^3 X_j X_j.
\label{invgen}
\eeq
\end{enumerate}

{\bf Step 2}: Let us calculate deformations of the potential 
$U(z,\bar z)$, the eigenfunction $\Psi(z,\bar z)$ and the constants
$C_j$ in formulas (\ref{weier}) corresponding to all infinitesimal 
conformal transformations of ${\Bbb R}^3$. 

\begin{enumerate}
\item Translations simply shift the constants $C_j$ and
change neither $U(z,\bar z)$ nor $\Psi(z,\bar z)$. Thus they
do not change $\Gamma[X,a,b]$, and without loss
of generality we can assume
\beq
C_j=0, \ j=1,2,3.
\eeq
\item A simple direct calculation based on thr representation (\ref{weier2}) 
(we do not like to reproduce it here) shows
that the rotations in ${\Bbb R}^3$ correspond to the following
transformations of the eigenfunction $\vec\Psi(z,\bar z)$
\beq
\left(\barr{c}\Psi_1(z,\bar z)\\ \Psi_2(z,\bar z)\earr\right)\rightarrow
\alpha \left(\barr{c}\Psi_1(z,\bar z)\\ \Psi_2(z,\bar z)\earr\right)+
\beta \left(\barr{c}\bar\Psi_2(z,\bar z)\\ -\bar\Psi_1(z,\bar z)\earr\right), 
|\alpha|^2+|\beta|^2=1
\label{rotations}
\eeq
where $\alpha$ and $\beta$ are some complex parameters.

Both functions $\vec\Psi(z,\bar z)$ and 
\beq
\vec\Psi^+(z,\bar z)=
\left(\barr{c}\bar\Psi_2(z,\bar z)\\ -\bar\Psi_1(z,\bar z)\earr\right), 
\label{conjpsi}
\eeq
satisfy the Dirac equation (\ref{diraceq}) with the same potential 
$U(z,\bar z)$. Thus the rotations does not change the potential 
and $\Gamma[X,a,b]$. 
\item The dilation is generated by the scaling transform 
\beq
\delta\vec\Psi(z,\bar z)= \frac12 \vec\Psi(z,\bar z)
\label{dilpsi}
\eeq
and changes neither $U(z,\bar z)$ nor $\Gamma[X,a,b]$.
\item The only nontrivial transformations of the Dirac operator correspond 
to the inversion generators. Up to conjugations by rotations all
generators (\ref{invgen}) are equivalent. Thus it is sufficient to
prove, that $\Gamma[X,a,b]$ is invariant if we apply the 
generator $-K_3$:
\beq
\delta X_1=-2 X_1 X_3,\ \delta K_2=-2 X_2 X_3,\ \delta X_3=
-X_3^2+X_1^2+X_2^2.
\label{maininv}
\eeq

Let us introduce the following notation
\beq
W(z,\bar z)=X_1(z,\bar z)-i X_2(z,\bar z).
\eeq

The transformation (\ref{maininv}) corresponds to
the following transformation of the function $\vec\Psi(z,\bar z)$
\beq
\barr{c}
\delta \Psi_1(z,\bar z)= -X_3(z,\bar z) \Psi_1(z,\bar z) + 
i W(z,\bar z) \bar\Psi_2(z,\bar z)
\\
\delta \Psi_2(z,\bar z)= -X_3(z,\bar z) \Psi_2(z,\bar z) - 
i W(z,\bar z) \bar\Psi_1(z,\bar z)
\earr
\label{invpsi}
\eeq

Let us check it.
$$
\barr{c}
\delta X_3(z,\bar z)=
-\int\limits_{z_0}^{z} \left[\left(\delta\Psi_2(z',\bar z')
\bar\Psi_1(z',\bar z') + \Psi_2(z',\bar z')\delta\bar\Psi_1(z',\bar z') 
\right)dz' \right. +
\\ 
+\left. \left( \delta\Psi_1(z',\bar z')\bar\Psi_2(z',\bar z')+
\Psi_1(z',\bar z')\delta\bar\Psi_2(z',\bar z') \right) d\bar z' \right]=
\earr
$$
$$
\barr{c}
=-\int\limits_{z_0}^{z} \left[\left(
-2 X_3(z',\bar z') \Psi_2(z',\bar z') \bar\Psi_1(z',\bar z') - \right. \right.
\\  
\left . -i W(z',\bar z') \bar\Psi_1^2(z',\bar z') -
i \bar W(z',\bar z') \Psi_2^2(z',\bar z')
\right)dz' +
\\
+ \left( 
-2 X_3(z',\bar z') \Psi_1(z',\bar z') \bar\Psi_2(z',\bar z') + \right.
\\
\left. \left. +i W(z',\bar z') \bar\Psi_2^2(z',\bar z') +
i \bar W(z',\bar z') \Psi_1^2(z',\bar z')
\right) d\bar z' \right]=
\earr
$$
$$
\barr{c}
=-\int\limits_{z_0}^{z} \left[\left( \vphantom{\bar \Psi}
2 X_3(z',\bar z') (\partial_{z'} X_3(z',\bar z')) - \right. \right.
\\
\left . - W(z',\bar z') (\partial_{z'}\bar W(z',\bar z') ) -
\bar W(z',\bar z') (\partial_{z'} W(z',\bar z'))
\right)dz' +
\\
+ \left( \vphantom{\bar \Psi}
2 X_3(z',\bar z') (\partial_{\bar z'} X_3(z',\bar z')) - \right.
\\
\left. \left. - W(z',\bar z') (\partial_{\bar z'}\bar W(z',\bar z')) -
\bar W(z',\bar z') (\partial_{\bar z'}W(z',\bar z'))
\right) d\bar z' \right]=
\earr
$$
$$
\barr{c}
=-\int\limits_{z_0}^{z} \left[\left( \vphantom{\bar \Psi}
\partial_{z'} X_3^2(z',\bar z')
-\partial_{z'} (W(z',\bar z') \bar W(z',\bar z') )
\right)dz' \right.+
\\
+ \left. \left( \vphantom{\bar \Psi}
\partial_{\bar z'} X_3^2(z',\bar z')-
\partial_{\bar z'}(W(z',\bar z')\bar W(z',\bar z'))
\right) d\bar z' \right]=
\earr
$$
$$
=W(z,\bar z)\bar W(z,\bar z)-X_3^2(z,\bar z)
$$
$$
\barr{c}
\delta W(z,\bar z)=
i\int\limits_{z_0}^{z}2 \left[\vphantom{\bar\Psi}
\delta\Psi_2(z',\bar z')\Psi_2(z',\bar z') dz'- 
\delta\Psi_1(z',\bar z') \Psi_1(z',\bar z')d\bar z' \right]=
\earr
$$
$$
\barr{c}
=\int\limits_{z_0}^{z}2 \left[\left(\vphantom{\bar\Psi}
-X_3(z',\bar z')i \Psi_2^2(z',\bar z') + 
W(z',\bar z') \bar\Psi_1(z',\bar z')\Psi_2(z',\bar z')\right) dz' \right.+
\\
+\left.\left(\vphantom{\bar\Psi}
X_3(z',\bar z')i \Psi_1^2(z',\bar z') + 
W(z',\bar z') \bar\Psi_2(z',\bar z')\Psi_1(z',\bar z') 
\right)d\bar z' \right]=
\earr
$$
$$
\barr{c}
=2\int\limits_{z_0}^{z} \left[\left(\vphantom{\bar\Psi}
-X_3(z',\bar z')(\partial_{z'} W(z',\bar z')) - 
W(z',\bar z') (\partial_{z'}X_3(z',\bar z'))\right) dz' \right.+
\\
+\left.\left(\vphantom{\bar\Psi}
- X_3(z',\bar z')(\partial_{\bar z'} W(z',\bar z'))
- W(z',\bar z')(\partial_{z'}X_3(z',\bar z'))  
\right)d\bar z' \right]=
\earr
$$
$$
=-2\int\limits_{z_0}^{z} \left[\left(\vphantom{\bar\Psi}
\partial_{z'} (X_3(z',\bar z')W(z',\bar z')) \right)d z' +
\left(\vphantom{\bar\Psi}
(\partial_{\bar z'}(X_3(z',\bar z') W(z',\bar z'))
\right)d\bar z' \right]=
$$
$$
=-2 W(z,\bar z)X_3(z,\bar z)
$$

The corresponding transformation of the potential $U(z,\bar z)$ reads as
\beq
\delta U(z,\bar z)=\Psi_1(z,\bar z) \bar \Psi_1(z,\bar z)-
\Psi_2(z,\bar z) \bar \Psi_2(z,\bar z).
\label{defpot}
\eeq
\end{enumerate}

{\bf Step 3}: Let us  calculate the the deformation of the
Bloch functions corresponding to (\ref{defpot}). 

Let $\vec\psi(\lambda,z,\bar z)$ be the Bloch function of $L$. For any
$\lambda$ such that at least one of the functions $\hbox{Im}\,p_x(\lambda)$, 
$\hbox{Im}\,p_y(\lambda)$ is not equal to zero (``non-physical'' $\lambda$)
define the following pair of functions
\beq
\barr{c}
\Omega_1(\lambda,z,\bar z) = \int\limits_{\infty}^z
\psi_2(\lambda,z',\bar z')\bar\Psi_1(z',\bar z')dz'+
\psi_1(\lambda,z',\bar z')\bar\Psi_2(z',\bar z')d\bar z'
\\
\Omega_2(\lambda,z,\bar z)= \int\limits_{\infty}^z
\psi_2(\lambda,z',\bar z')\Psi_2(z',\bar z')dz'-
\psi_1(\lambda,z',\bar z')\Psi_1(z',\bar z')d\bar z'.
\earr
\label{omegas}
\eeq
where the integrals are taken along an arbitrary path in the $z$-plane, 
connecting the points $z$ and $\infty$ such that the integrand decays 
exponentially along this path. The integrands in (\ref{omegas})
are closed 1-forms thus the integrals do not depend on a concrete
choice of the path. Using the same arguments as in Section
\ref{sec:CBA} below we may easily prove that $\Omega_1(\lambda,z,\bar{z})$
and $\Omega_2(\lambda,z,\bar{z})$ are meromorphic in 
$\lambda$ outside infinity.

The function $\vec\Psi(z,\bar z)$ is double-periodic or anti-periodic
in $z$ (see (\ref{persol})), thus
the functions $\Omega_1(\lambda,z,\bar z)$, $\Omega_2(\lambda,z,\bar z)$
have the following periodicity properties: 
\beq
\barr{c}
\Omega_k(\lambda,z+1,\bar z+1)=
{\cal W}_1 w_1(\lambda)\Omega_k(\lambda,z,\bar z)
\\ 
\Omega_k(\lambda,z+\tau,\bar z+\bar\tau)=
{\cal W}_2 w_2(\lambda)\Omega_k(\lambda,z,\bar z)
\earr
\ \ k=1,2.
\label{omegaper}
\eeq 
(see (\ref{persol}) for the definition of ${\cal W}_1$, ${\cal W}_2$.)
\begin{lemma}
\label{defgamma}
The variation of the Bloch function $\vec\psi(\lambda,z,\bar z)$ 
corresponding to (\ref{defpot}) reads as
\beq
\barr{c}
\delta \psi_1(\lambda,z,\bar z)= \Omega_1(\lambda,z,\bar z) \Psi_1(z,\bar z) - 
\Omega_2(\lambda,z,\bar z) \bar\Psi_2(z,\bar z)+ 
\alpha(\lambda)\psi_1(\lambda,z,\bar z)
\\
\delta \psi_2(\lambda,z,\bar z)= \Omega_1(\lambda,z,\bar z) \Psi_2(z,\bar z) +
\Omega_2(\lambda,z,\bar z) \bar\Psi_1(z,\bar z) + 
\alpha(\lambda)\psi_2(\lambda,z,\bar z)
\earr
\label{defbloch}
\eeq
where $\alpha(\lambda)$ is some meromorphic function, fixed by the 
normalization condition:
\beq
\left.\delta\psi_1(\lambda,z,\bar z)+\delta\psi_2(\lambda,z,\bar z)\right|_
{z=z_1}=0,
\label{norm_1}
\eeq
\end{lemma}

Proof of Lemma \ref{defgamma}. To start with, let us recall a simple
fact from Bloch theory (see for example \cite {Kri89}).

Generically, if we calculate variations of the Bloch function, we
deform $\Gamma$, and we can not assume both $\delta w_1(\lambda)=0$ and
$\delta{w_2(\lambda)}=0$ simultaneously. To compare functions on
different subvarieties
in ${\Bbb C}^2$ we have to fix some connection. The simplest way to do
this is to assume $\delta w_1(\lambda)=0$.

Then the variation of the Bloch function can be found as the unique 
solution of the linearized Dirac equation
\beq
\delta L \vec\psi(\lambda,z,\bar{z})+L \delta\vec\psi(\lambda,z,\bar{z})=0
\label{linearized_dirac}
\eeq
satisfying (\ref{norm_1}) such that
\beq
\delta\vec\psi(\lambda,t_1+\tau t_2,t_1+\bar\tau t_2)=
O\left(\left[1+\left|t_2\right|\right]e^{ip_1(\lambda)t_1+ip_2(\lambda)t_2}
\right).
\label{delta_psi_estimate}
\eeq

A simple direct calculation shows that (\ref{defbloch}) solves 
(\ref{linearized_dirac}). From (\ref{omegaper}) and (\ref{persol}) it
follows, that
\beq
\barr{c}
\delta\vec\psi(w_1,w_2,z+1,\bar z+1)=w_1\delta\vec\psi(w_1,w_2,z,\bar z),
\\ 
\delta\vec\psi(w_1,w_2,z+\tau,\bar z+\bar\tau)=
w_2\delta\vec\psi(w_1,w_2,z,\bar z),
\earr
\label{per_delta_bloch}
\eeq
thus variations of the type (\ref{defbloch}) satisfy 
(\ref{delta_psi_estimate}). This completes the proof.

The function $\delta\psi(\lambda,z,\bar{z})$ defined by
(\ref{defbloch}) has the same periodicity properties as the original
Bloch function $\psi(\lambda,z,\bar{z})$ (see formulas
(\ref{blochsol}) and (\ref{per_delta_bloch}) respectively). Thus our
special variations satisfy $\delta w_1(\lambda)=0$ and  
$\delta w_2(\lambda)=0$ simultaneously.

{\bf Step 4}: From $(\ref{defbloch})$ it follows, that if we apply 
infinitesimal  conformal transformation $-K_3$ to our immersion, 
the Bloch functions of the deformed Dirac operator are meromorphic on
the same variety $\Gamma[X,a,b]$ as the original Bloch functions and have the
same multipliers $w_1$, $w_2$. Thus our deformation does not change
$\Gamma[X,a,b]$. This completes the proof of Theorem~\ref{main_theorem}.

\begin{example}  {\bf Surfaces of revolution}. Let $\gamma$ be a
  closed non-selfintersecting curve in the half-plane $X_2=0$, $X_1>0$
  in ${\Bbb R}^3$. Rotating $\gamma$ about the axes $X_1=X_2=0$ we get
  a surface of revolution. It is always a torus with a fixed pair of
  periods. 

  Such surfaces are essentially simpler from the soliton point of
  view. Potentials $U(z,\bar{z})$ corresponding to such surfaces
  depend only on one real variable $x=\,\hbox{Re}\,z$. Instead of the
  fixed energy spectral transform for the 2-dimensional Dirac operator
  we have the spectral transform for the $2\times2$ first-order matrix
  differential operator in one variable. Periodic direct spectral
  transform for such operators (for both finite-gap and infinite-gap
  potentials) was developed by one of the authors (M.S.) in
  \cite{Sch95}. 

  The MNV equations for surfaces of revolution are reduced to the
  well-studied Modified Korteweg-de~Vries equations (MKdV). In
  contrast with MNV all higher MKdV conservation laws are local in
  terms of the potential. Equation (\ref{defpot}) in this situation
  was integrated by V.~K.~Melnikov in \cite{Mel90} in the class of
  potentials sufficiently fast decaying at infinity. Our theorem for 
  the surfaces of revolution does not follow formally from
  \cite{Mel90} because the periodic MKdV theory and the decaying at
  infinity one use different technical tools. Nevertheless it is 
  possible to simplify essentially our proof in this specific case.
\end{example}

\label{sec:cons}
\section{Appendix. Modified Novikov-Veselov equations with periodic
  boundary conditions.}

In this Appendix we discuss the zero-energy spectral transform for
the double-periodic Dirac operator (and the Generalized Weierstrass transform
respectively) from the soliton point of view. This transform is
naturally connected with a completely integrable hierarchy of
integro-differential equations with 2 spatial variables known as Modified
Novikov-Veselov hierarchy (MNV).

Formally, the results of this Section are not used in our proof of 
Theorem~\ref{main_theorem}, but we hope they allow the reader to
gain a better understanding of the problem. 

There is a rather big amount of papers dedicated to the periodic
problem for soliton equations (see for example the textbook
\cite{ZMNP}). Nevertheless, the direct spectral transform for 
two-dimensional Dirac operator was never studied in such context in
the literature known to us. Important properties of Bloch varieties for 
multidimensional Dirac operators were proved in \cite{Kuch93}, but 
they are not sufficient for the purpos of integrating the periodic
MNV equations. 

From the point of view of the Bloch theory the two-dimensional 
double-periodic Dirac  operator is rather similar to 
the two-dimensional double-periodic Schr\"odinger operator. 
The fixed-energy direct spectral transform for last one was 
constructed by soliton methods by I.~M.~Krichever \cite{Kri89}.
This problem was studied more detail by J.~Feldman, H.~Kn\"orrer and 
E.~Trubowitz (see \cite{FKT}). In Section~\ref{sec:asymtotic} 
we describe the structure of the Dirac Bloch variety by methods
analogous to \cite{Kri89}. It is important to remark, that in spite of
the similarity between these two problems we have to overcome some
additional technical difficulties on this way. The asymptotic
expansion of the Bloch variety gives us ``quasi-local'' MNV
conservation laws. 
 
The zero-energy scattering problem for the two-dimensional Dirac
operator possesses an infinite-dimensional algebra of symmetries,
generated by the MNV equations. They were constructed by L.~V.~Bogdanov 
in \cite{Bog87}. In \cite{Bog87} a generalization of the Miura transform 
was defined, and it was shown, that this transform maps the MNV
equations to the Novikov-Veselov hierarchy (see \cite{VN84}),
associated with the fixed-energy two-dimensional Schr\"odinger
operator. MNV equations in the space of functions, decaying at
infinity were integrated by the so-called method of
$\bar\partial$-problem in \cite{Bog87}. Periodic MNV theory is
discussed in Section~\ref{sec:mnv}.

In contrast with the one-dimensional soliton systems, the
two-dimensional ones essentially depend on the boundary conditions. To
define the periodic MNV hierarchy uniquely we have to fix some
constants of integration. One of the simplest way to do it is to define
the MNV hierarchy in terms of the so-called Cauchy-Baker-Akhiezer (CBA)
kernel. This kernel was introduced for the Kadomtsev-Petviashvily
hierarchy by A.~Yu.Orlov and one of the authors (P.G.) in \cite{GO89}
(see also \cite{GOS93}). MNV hierarchy in terms of the CBA kernel is
discussed in Section~\ref{sec:CBA}.

In this Section we always assume that $U(z,\bar z)$ is real, smooth
and double-periodic (\ref{u}).

\subsection{Asymptotic structure of the Bloch variety.}
\label{sec:asymtotic}

For large $\hbox{Im}\,p_1$, $\hbox{Im}\,p_2$ the structure of the
Bloch variety can be studied by methods of perturbation
theory. Following \cite{Kri89} we start from the Dirac operator with
zero potential $U(z,\bar z)\equiv0$. The corresponding Bloch variety
is a union of two Riemann spheres
$\Gamma^{(0)}=\Gamma_1^{(0)}\cup\Gamma_2^{(0)}$, 
$\Gamma_1^{(0)}=\Gamma_2^{(0)}={\Bbb {CP}}^1$ with a coordinate $\lambda$
and
\beq
\vec\psi(\lambda,z,\bar z)=\left(\barr{c} 1 \\ 0 \earr \right )
e^{\lambda \bar z}, \ \lambda\in\Gamma_1^{(0)}, \
\vec\psi(\lambda,z,\bar z)=\left(\barr{c} 0 \\ 1 \earr \right )
e^{\lambda z}, \ \lambda\in\Gamma_2^{(0)}.
\label{zeropotential}
\eeq

A pair $\lambda_1\in\Gamma_1^{(0)}$, $\lambda_2\in\Gamma_2^{(0)}$ is
called resonant if 
\beq
e^{\lambda_1-\lambda_2}=1,\ e^{\lambda_1\bar\tau-\lambda_2\tau}=1,
\label{resonant}
\eeq
and non-resonant otherwise. All resonant pairs are given by the
following formulas:
\beq
\lambda_1^{(m,n)}=\frac{\pi m \hbox{Re}\,\tau -\pi n}
{\hbox{Im}\,\tau}+i\pi m,\ 
\lambda_2^{(m,n)}=\bar\lambda_1^{(m,n)},\ m,n\in\Bbb Z.
\label{respairs}
\eeq
Let us call a point $\lambda\in\Gamma_1^{(0)}$ non-resonant if the pair
$\lambda\in\Gamma_1^{(0)}$ and $\bar\lambda\in\Gamma_2^{(0)}$ is
non-resonant. The antiholomorphic involution $\sigma\theta$ maps 
$\Gamma_2^{(0)}$ to $\Gamma_1^{(0)}$ thus it is sufficient to
develop a perturbation theory only on $\Gamma_1^{(0)}$. 

Let $\varepsilon$, $R$ be some positive constants. Denote by 
$\Gamma^{(0)}_{\varepsilon,R}$ the domain, obtained from
$\Bbb{CP}^1$ by removing $\varepsilon$ neighbourhoods of all
resonant points $\lambda_1^{(m,n)}$ and the disk $|\lambda|\le R$.

\begin{lemma}
\label{lemma:asymptotic}
For any $\varepsilon >0$ there exists a constant $R(\varepsilon)$
such that in the domain $\Gamma^{(0)}_{\varepsilon,R(\varepsilon)}$
there exists an unique solution of the Dirac equation (\ref{diraceq})
with normalization (\ref{norm}) such that
\beq
\barr{c}
\vec\psi(\lambda,z+1,\bar z+1)=e^{\lambda+h(\lambda)}
\vec\psi(\lambda,z,\bar z),
\\ 
\vec\psi(\lambda,z+\tau,\bar z+\bar\tau)=e^{\lambda\bar\tau+h(\lambda)\tau}
\vec\psi(\lambda,z,\bar z),
\earr
\label{pertbloch}
\eeq
where $h(\lambda)$ is uniquely defined under the condition that
$h(\lambda)\rightarrow 0$ as $\lambda\rightarrow\infty$. The functions
$\vec\psi(\lambda,z,\bar z)$, and $h(\lambda)$ are holomorphic in
$\lambda$ in the domain $\Gamma^{(0)}_{\varepsilon,R(\varepsilon)}$.
\end{lemma}
A proof of this statement is to appear in a forthcomming paper by
I.~A.~Taimanov and one of the authors (M.S.). It is rather long and
rather technical. We do not want to present it in our text.
 
\begin{statement}
\label{st2}
The functions $\vec\psi(\lambda,z,\bar z)$, $h(\lambda)$ defined in
Lemma~\ref{lemma:asymptotic} 
posses the following asymptotic expansions as $\lambda\rightarrow\infty$
\beq
\vec\psi(\lambda,z,\bar z)=e^{\lambda(\bar z-\bar z_1)+h(\lambda)(z-z_1)}
\left(\barr{c} 1+\frac{\phi_1(z,\bar z)}{\lambda}+
\frac{\phi_2(z,\bar z)}{\lambda^2}+
\frac{\phi_3(z,\bar z)}{\lambda^3}+
\frac{\phi_4(z,\bar z)}{\lambda^4}+\ldots
\\
\frac{\chi_1(z,\bar z)}{\lambda}+
\frac{\chi_2(z,\bar z)}{\lambda^2}+
\frac{\chi_3(z,\bar z)}{\lambda^3}+
\frac{\chi_4(z,\bar z)}{\lambda^4}+\ldots
\earr\right)
\label{phias}
\eeq
\beq
h(\lambda)=\frac{h_1}{\lambda}+\frac{h_3}{\lambda^3}+\frac{h_5}{\lambda^5}+
\ldots
\label{has}
\eeq
(the Bloch variety $\Gamma$ has a symmetry 
$\sigma: (w_1,w_2)\rightarrow(w_1^{-1},w_2^{-1})$, 
$\sigma(\lambda)=-\lambda$, thus   
all even coefficients in (\ref{has}) are identically zero.)
\end{statement}
Unfortunately, at this moment we do not have any complete proof of
this Statement. It is rather clear how to do it, but this proof needs
a rather big piece of asymptotic analysis, and we are not in a position
to do it presently. But we know, that it is fulfilled at least in two
important specific situations:
\begin{enumerate}
\item $U(z,\bar{z})$ is an algebraic-geometrical (or, equivalently,
  finite-gap) potential. It means, that the normalization of the 
  zero-energy Bloch variety is algebraic (has finite genus).
\item $U(z,\bar{z})$ depend only on one real variable
  $x=\,\hbox{Re}\,z$. This fact was proved by one of the authors
  (M.S.) in \cite{Sch95}. Such potentials corresponds to surfaces of
  revolution. 
\end{enumerate}
We would like to remark, that in \cite{FKT} a class of Riemann surfaces was
introduced, which are in some sense similar to compact Riemann surfaces,
and the zero-energy level of a two-dimensional Schr\"odinger operator
belongs to this class. Therefore it is natural to expect that the
zero-energy level of the two-dimensional Dirac operator also belongs
to this class.

To define the Modified Novikov-Veselov equations (MNV)
and their conservation laws it is sufficient to have a formal solution 
of the Dirac equation in the form (\ref{phias}), (\ref{has}). Let us
show that such solution exists and is unique if we assume that all 
$\phi_k(z,\bar z)$, $\chi_k(z,\bar z)$ are bounded in the whole
$z$-plane.

Inserting (\ref{phias}) and 
\beq
h(\lambda)=\frac{h_1}{\lambda}+\frac{h_2}{\lambda^2}+\frac{h_3}{\lambda^3}+
\ldots
\label{has2}
\eeq
in (\ref{diraceq}) we get the following system of equations
\beq
\barr{c}
\chi_1(z,\bar z)=-U(z,\bar z),
\\
\chi_k(z,\bar z)=-\partial_{\bar z}\chi_{k-1}(z,\bar z)-U(z,\bar z)
\phi_{k-1}(z,\bar z),\ k>1,
\\
\partial_z\phi_k(z,\bar z)=U(z,\bar z)\chi_{k}(z,\bar z)
-h_k-\sum_{j=1}^{k-1}h_j\phi_{k-j}(z,\bar z).
\earr
\label{recceq}
\eeq
We solve this system by induction. First we find $\chi_1(z,\bar z)$,
then $\phi_1(z,\bar z)$, then $\chi_2(z,\bar z)$, then 
$\phi_2(z,\bar z)$ and so on. To find $\chi_k(z,\bar z)$ at each step
we differentiate some double-periodic functions, obtained at previous
steps. Thus they are defined uniquely and are automatically
double-periodic. To find $\phi_k(z,\bar z)$ we have to invert the 
operator $\partial_z$ in the space of functions, bounded in the whole 
$z$-plane (any bounded solution $\phi_k(z,\bar z)$ is automatically 
double-periodic). This is possible if and only if the mean value of the
right-hand side is equal to zero:
\beq
<<U(z,\bar z)\chi_{k}(z,\bar z)-h_k
-\sum_{j=1}^{k-1}h_j\phi_{k-j}(z,\bar z)>>=0.
\label{solv}
\eeq
where:
\beq
<<F(z,\bar z)>>=\int_0^1\int_0^1dt_1dt_2 F(t_1+\tau t_2,t_1+\bar\tau t_2).
\eeq
Thus at each step we find $h_k$ from (\ref{solv}) and then 
calculate $\phi_k(z,\bar z)$.
The function $\phi_k(z,\bar z)$ is determined by the system
(\ref{recceq}) uniquely up to adding an arbitrary constant. 
This constant is fixed by the normalization condition (\ref{norm}). 

Let us check that the function $h(\lambda)$ does not depend on
the normalization point $z_1$. If we change the point $z_1$, then
we change the integration constants. But these constants can be 
arbitrary shifted by multiplying the whole solution to a formal series
in $\lambda$
\beq
\left[\!\barr{c} 1+\frac{\phi_1(z,\bar z)}{\lambda}+
\frac{\phi_2(z,\bar z)}{\lambda^2}+\ldots
\\
\frac{\chi_1(z,\bar z)}{\lambda}+
\frac{\chi_2(z,\bar z)}{\lambda^2}+\ldots
\earr\!\right]
\!\rightarrow\!
\left[1\!+\!\frac{\alpha_1}{\lambda}\!+\!\frac{\alpha_2}{\lambda^2}\!+\!
\ldots \right]
\left[\!\barr{c} 1+\frac{\phi_1(z,\bar z)}{\lambda}+
\frac{\phi_2(z,\bar z)}{\lambda^2}+\ldots
\\
\frac{\chi_1(z,\bar z)}{\lambda}+
\frac{\chi_2(z,\bar z)}{\lambda^2}+\ldots
\earr\!\right].
\label{changenorm}
\eeq
This multiplication does not affect $h(\lambda)$.

We have proved, that the constants $h_1$, $h_3$, $h_5$, \ldots are
completely determined by the potential $U(z,\bar z)$ and all
$h_{2k}=0$. Thus we have constructed an infinite sequence of functionals
$h_{2k+1}[U]$. The formulas for the first two of them are:
\beq
h_1=-<<U^2(z,\bar z)>>
\label{h_1}
\eeq
\beq
h_3=-<<U(z,\bar z)U_{\bar z\bar z}(z,\bar z)-(U^2(z,\bar z)+h_1)
V_{1\bar z}(z,\bar z)>>,
\label{h_3}
\eeq
where
\beq
V_1(z,\bar z)=\partial^{-1}_z (U^2(z,\bar z)+h_1).
\label{V_1}
\eeq
Let us point out that adding an arbitrary constant to $V_1(z,\bar z)$
does not affect $h_3$. 

We have defined an infinite collection of functionals $h_{2k+1}[U]$, 
$k=0,1,2,\ldots$. The following Statement explains why these functionals 
are so important.
\begin{statement}
\label{lemma:conservation}
The quantities $h_{2k+1}[U]$ are {\bf conservation laws for the whole
hierarchy of the Modified Novikov-Veselov equations (MNV)}.
\end{statement}
Assuming that Statement~\ref{st2} is fulfilled we the prove this
Statement at the end of the Section.

It is well-known in the soliton theory that integrable systems with 
one spatial variable usually have infinitely many local conservation 
laws. For multidimensional soliton systems we normally have 
the opposite situation: almost all conservation laws are nonlocal. 
Let us briefly discuss the case of the MNV hierarchy.

A functional $Q[U]$ is called local if it possesses the following 
representation:
\beq
Q[U]=<<q(U,U_z,U_{\bar z},U_{z z},U_{z\bar z},U_{\bar z\bar z},\ldots)>>
\label{local_functionals}
\eeq
where the density $q(\ldots)$ depend only on $U(z,\bar z)$ and finite 
number of its derivatives. Of course $h_1[U]$ is local. The next 
conservation law $h_3[U]$ is nonlocal  
because the corresponding density depends on an auxiliary function 
$V_{1\bar z}(z,\bar z)$, and to calculate $V_{1\bar z}(z,\bar z)$
we have to know $U(z,\bar z)$ on the whole $z$-plane. It it rather 
evident that higher functionals $h_{2k+1}[U]$ are also nonlocal. 
This nonlocality creates no serious problems if the potential is
double-periodic, but it is very difficult to extend existing definitions
to wider classes of boundary conditions.

In \cite{Kri89} perturbation theory was developed also in 
neighbourhoods of resonant pairs. It can be shown that for sufficiently
large $\lambda$ the surface $\Gamma$ is obtained from the $\Gamma{(0)}$
by attaching small handles to the resonant pairs. We do not use this  
fact, thus we do not want to discuss it now. However we use the following
property of $\Gamma$:

\begin{corollary} The Bloch variety $\Gamma$ has two
infinite points, corresponding to the points $\lambda=\infty$ in 
$\Gamma_1^{(0)}$ and $\Gamma_2^{(0)}$. We denote them by $\infty_+$ and
$\infty_-$ respectively. 
\end{corollary}
This statement follows immediately from Lemma \ref{lemma:asymptotic}.

Further we shall use the following notation:
$\lambda\rightarrow\infty_+$, where $\lambda$ is a point of $\Gamma$, 
always means that $\lambda$ tends to $\infty_+$ in the domain
$\Gamma^{(0)}_{\varepsilon,R(\varepsilon)}$; the notation
$\lambda\rightarrow\infty_-$ always means that
$\bar\lambda\in\Gamma^{(0)}_{\varepsilon,R(\varepsilon)}$.

\subsection{Cauchy-Baker-Akhiezer kernel.}
\label{sec:CBA}

To study symmetries of the soliton equation it is convenient to use the 
Cauchy-Baker-Akhiezer kernel (CBA) (see  \cite{GO89}).
In this Section we define the CBA kernel on the zero-energy
Bloch variety of the two-dimensional Dirac operator.

To start with, suppose $(\lambda,\mu)$ is a pair of points in $\Gamma$
such that
\begin{enumerate}
\item $\vec\psi(\nu,z,\bar z)$ is non-singular at the points
$\nu=\lambda$ and $\nu=\sigma\mu$. (Let us recall, that outside
$\nu=\infty$ the poles of $\vec\psi(\nu,z,\bar z)$ arose due to the 
normalization (\ref{norm}) and do not depend on $z$ and $\bar z$.)
\item At least one of the following conditions is fulfilled:
\beq
\hbox{Im}\,p_1(\lambda)-\hbox{Im}\,p_1(\mu)\ne 0 \ \ \hbox{or} \ \ 
\hbox{Im}\,p_2(\lambda)-\hbox{Im}\,p_2(\mu)\ne 0.
\label{non-resonant_1}
\eeq
\end{enumerate} 
Then we may define $\tilde\omega(\lambda,\mu,z,\bar z)$ by: 
\beq
\tilde\omega(\lambda,\mu,z,\bar z)=\int_{\infty}^z 
d\tilde\omega(\lambda,\mu,z',\bar z'), \ 
\label{cba_1}
\eeq
where
\beq
d\tilde\omega(\lambda,\mu,z',\bar z')=
\psi_2(\lambda,z',\bar z')\psi_2(\sigma\mu,z',\bar z')dz'-
\psi_1(\lambda,z',\bar z')\psi_1(\sigma\mu,z',\bar z')d\bar z'
\label{cba_2}
\eeq
The integral in (\ref{cba_1}) is taken along some path $\gamma$ in the
$z$-plane such that
\begin{enumerate}
\item $\gamma$ connects the point $z$ with $\infty$.
\item The form $d\tilde\omega(\lambda,\mu,z',\bar z')$ decays
exponentially as $z'\rightarrow\infty$ along $\gamma$. 
\end{enumerate}
Condition (\ref{non-resonant_1}) guarantees the existence of such path.

From the Dirac equation (\ref{diraceq}) if follows that the form 
$d\tilde\omega(\lambda,\mu,z',\bar z')$ is closed that the integral 
(\ref{cba_1}) is well-defined.

The next step is to show, that for a fixed $\mu$ our function 
$\tilde\omega(\lambda,\mu,z,\bar z)$ is meromorphic in $\lambda$ on $\Gamma$
outside the points $\infty_+$ and $\infty_-$. 

Suppose, that a pair $(\lambda,\mu)$ satisfies the following
condition, which is much weaker than (\ref{non-resonant_1}).

\medskip
$2'$. At least one of the following combinations 
\beq
\frac1{2\pi}\left(p_1(\lambda)-p_1(\mu)\right) \ \hbox{or} \ 
\frac1{2\pi|\tau|}\left(p_2(\lambda)-p_2(\mu)\right)
\ \hbox{is non-integer,}
\label{non-resonant_2}
\eeq
or, equivalently, at least one of the following inequalities is
fulfilled
\beq
w_1(\lambda)w_1^{-1}(\mu)\ne 1, \ \hbox{or} \ 
w_2(\lambda)w_2^{-1}(\mu)\ne 1.
\label{non-resonant_3}
\eeq
Assume that $w_1(\lambda)w_1^{-1}(\mu)\ne 1$. Then we
have either 
\beq
\tilde\omega(\lambda,\mu,z,\bar z)=\int_{-\infty}^{0} 
\tilde\omega_x(\lambda,\mu,z+t,\bar z+t)dt \
\label{cba_x+}
\eeq
or
\beq
\tilde\omega(\lambda,\mu,z,\bar z)=-\int^{\infty}_{0} 
\tilde\omega_x(\lambda,\mu,z+t,\bar z+t)dt,
\label{cba_x-}
\eeq
where
\beq
\tilde\omega_x(\lambda,\mu,z',\bar z')=
\psi_2(\lambda,z',\bar z')\psi_2(\sigma\mu,z',\bar z')-
\psi_1(\lambda,z',\bar z')\psi_1(\sigma\mu,z',\bar z'),
\label{cba_omega_x}
\eeq
and $t\in{\Bbb R}$.

If $\left|w_1(\lambda)w_1^{-1}(\mu)\right|\ge 1$ from (\ref{cba_x+}) we get:
\beq
\barr{c}
\tilde\omega(\lambda,\mu,z,\bar z)=
\\
=\left(\int\limits_{-1}^{0}+
\int\limits_{-2}^{-1}+
\int\limits_{-3}^{-2}+
\ldots \right)
\tilde\omega_x(\lambda,\mu,z+t,\bar z+t)dt =
\\
=\left(\frac{w_1(\mu)}{w_1(\lambda)}+\frac{w_1^2(\mu)}{w_1^2(\lambda)}+
\frac{w_1^3(\mu)}{w_1^3(\lambda)}+\ldots  \right)
\int\limits_{0}^{1}
\tilde\omega_x(\lambda,\mu,z+t,\bar z+t)dt =
\\
=\frac{w_1(\mu)}{w_1(\lambda)-w_1(\mu)}
\int\limits_{0}^{1}
\tilde\omega_x(\lambda,\mu,z+t,\bar z+t)dt
\earr
\label{cba_reg}
\eeq
Similarly if $\left|w_1(\lambda)w_1^{-1}(\mu)\right|\le 1$ 
from (\ref{cba_x-}) we get:
\beq
\barr{c}
\tilde\omega(\lambda,\mu,z,\bar z)=
\frac{w_1(\mu)}{w_1(\lambda)-w_1(\mu)}
\int\limits_{0}^{1}
\tilde\omega_x(\lambda,\mu,z+t,\bar z+t)dt
\earr
\label{cba_reg_2}
\eeq
The formulas (\ref{cba_reg}) and (\ref{cba_reg_2}) define meromorphic
continuations of the function $\tilde\omega(\lambda,\mu,z,\bar z)$
to the whole surface $\Gamma$. To conclude the proof, it remains to
note that formulas (\ref{cba_reg}) and (\ref{cba_reg_2}) coincide.

If $w_2(\lambda)w_2^{-1}(\mu)\ne 1$ then the integration
path can be chosen along the line $z'=z+\tau t$, and we obtain 
a new representation for the same function
$\tilde\omega(\lambda,\mu,z,\bar z)$:
\beq
\barr{c}
\tilde\omega(\lambda,\mu,z,\bar z)=
\frac{w_2(\mu)}{w_2(\lambda)-w_2(\mu)}
\int\limits_{0}^{1}
\tilde\omega_2(\lambda,\mu,z+\tau t,\bar z+\bar\tau t)dt,
\earr
\label{cba_reg_3}
\eeq
where
\beq
\tilde\omega_2(\lambda,\mu,z',\bar z')=
\tau\psi_2(\lambda,z',\bar z')\psi_2(\sigma\mu,z',\bar z')-
\bar\tau\psi_1(\lambda,z',\bar z')\psi_1(\sigma\mu,z',\bar z').
\label{cba_omega_2}
\eeq

Now we are in position to define the {\bf Cauchy-Baker-Akhiezer
kernel} $\omega(\lambda,\mu,z,\bar z)$:
\beq
\omega(\lambda,\mu,z,\bar z)=-\frac{1}{2\pi}\frac{dp_1(\mu)}
{<\tilde\omega_x(\mu,\mu,z,\bar z)>_x}
\tilde\omega(\lambda,\mu,z,\bar z)
\label{cba}
\eeq
where
\beq
<\tilde\omega_x(\mu,\mu,z,\bar z)>_x=
\int_{0}^{1}
\tilde\omega_x(\lambda,\mu,z+t,\bar z+t)dt.
\eeq
It is a simple exercise to check that the function
$<\tilde\omega_x(\mu,\mu,z,\bar z)>_x$ does not depend on $z$ and is
even in $\mu$, i.e.
\beq
<\tilde\omega_x(\mu,\mu,z,\bar z)>_x=
<\tilde\omega_x(\sigma\mu,\sigma\mu,z,\bar z)>_x.
\label{even}
\eeq
\beq
<\tilde\omega_x(\mu,\mu,z,\bar z)>_x=\mp1+\frac{c_2^{(\pm)}}{\mu^2}+
\frac{c_4^{(\pm)}}{\mu^4}+\ldots, \ \mu\rightarrow\infty_\pm
\label{even_2}
\eeq

\begin{lemma}
\label{lemma:cba}
The Cauchy-Baker-Akhiezer kernel $\omega(\lambda,\mu,z,\bar z)$
defined above has the following properties: 
\begin{enumerate}
\item For any fixed $z$, $\omega(\lambda,\mu,z,\bar z)$ is a
meromorphic function of $\lambda$ and a meromorphic 1-form in $\mu$
(both on the finite part of $\Gamma$).
\item For generic $\mu$, $\omega(\lambda,\mu,z,\bar z)$ has poles at
the poles of $\vec\psi(\lambda,z,\bar z)$ and at the point
$\lambda=\mu$. It is holomorphic outside these points at the finite part of
$\Gamma$. 
\item For a generic $\lambda$ and $\mu\rightarrow\lambda$
\beq
\omega(\lambda,\mu,z,\bar z)=\frac{1}{2\pi i}\frac{d\mu}{\mu-\lambda}+
\hbox{regular terms}, 
\label{residue_1}
\eeq
or, equivalently,
\beq
\oint_{\mu\in S} \omega(\lambda,\mu,z,\bar z)=1,
\label{residue_2}
\eeq
where $S$ is a small contour, surrounding the point $\lambda$.
\item Let $\mu\rightarrow\infty_+$  and let $\lambda$ be a fixed point
in $\Gamma$. Then we have the following formal expansion:
\beq
\omega(\lambda,\mu, z,\bar z)=\frac{d\mu}{2\pi i}
\left[\frac {\psi_1(\lambda,z,\bar z)}{\mu} + 
\sum_{k=2}^\infty \frac{R^{(+)}_{k}[\partial_{\bar z}]
\psi_1(\lambda,z,\bar z)}{\mu^k} \right]e^{-\mu\bar z-h(\mu)z},
\label{as_omega_+}
\eeq
where
\beq
R^{(+)}_{k}[\partial_{\bar z}]=
\partial_{\bar z}^{k-1}+
\sum_{l=0}^{k-2}v^{(+)}_{kl}(z,\bar z)\partial_{\bar z}^{l}
\label{r+}
\eeq 
are differential operators in $\bar z$ of order $k-1$. 
The coefficients $v^{(+)}_{kl}(z,\bar z)$ do not depend on $\lambda$ and
$\mu$, are differential polynomials of the asymptotic expansion 
coefficients  $\phi_j(z,\bar z)$, $\chi_j(z,\bar z)$, $j<k$. 
Each function $v^{(+)}_{kl}(z,\bar z)$ depends also on a finite number
of constants $h_{2j+1}$, $c^{(\pm)}_{2j}$.

Similarly for $\mu\rightarrow\infty_-$,  we have:
\beq
\omega(\lambda,\mu, z,\bar z)=
\frac{d\mu}{2\pi i}\left[\frac {\psi_2(\lambda,z,\bar z)}{\mu} + 
\sum_{k=2}^\infty \frac{R^{(-)}_{k}[\partial_{z}]
\psi_2(\lambda,z,\bar z)}{\mu^k}\right]e^{-\mu z-h(\mu)\bar z}.
\label{as_omega_-}
\eeq
where 
\beq
R^{(-)}_{k}[\partial_{z}]=
\partial_{z}^{k-1}+
\sum_{l=0}^{k-2}v^{(-)}_{kl}(z,\bar z)\partial_{z}^{l}.
\label{r-}
\eeq 
\item 
\beq
\barr{c}
\omega(\lambda,\mu,z+1,\bar z+1)=
w_1(\lambda)w_1^{-1}(\mu)\omega(\lambda,\mu,z,\bar z),
\\
\omega(\lambda,\mu,z+\tau,\bar z+\bar\tau)=
w_2(\lambda)w_2^{-1}(\mu)\omega(\lambda,\mu,z,\bar z).
\earr
\label{periods_omega}
\eeq
\end{enumerate}
\end{lemma}
\begin{statement}
\label{statement:cba}
If Statement~\ref{st2} is fulfilled, then the CBA kernel has the following
additional properties:
\begin{enumerate}
\item The formal expansions (\ref{as_omega_+}), (\ref{as_omega_-}) are
  asymptotic. 
\item \label{item:residues} Let $f^{(+)}(\lambda,z,\bar z)$ be the 
following formal series in $\lambda$
\beq
f^{(+)}(\lambda,z,\bar z)=\left\{
  \sum_{k=-N}^{\infty}\frac{f^{(+)}_k(z,\bar
    z)}{\lambda^k}\right\}e^{\lambda\bar z}.
\label{formal_+}
\eeq
Then
\beq 
\left.2\pi i\,\hbox{\rm res}\vphantom{\int}\right|_{\mu=\infty_+}
\!\!\!\!\!\!\!\!\omega(\lambda,\mu,z,\bar z)f^{(+)}(\mu,z,\bar z)= 
\left\{
  \sum_{k=-N}^{0}\frac{f^{(+)}_k(z,\bar
    z)}{\lambda^k}+O\left(\frac{1}{\lambda}\right)\right\}e^{\lambda\bar z}, 
\label{residue_++}
\eeq
as $\lambda\rightarrow\infty_+$ and 
\beq
\left.2\pi i\,\hbox{\rm res}\vphantom{\int}\right|_{\mu=\infty_+}
\!\!\!\!\!\!\!\!\omega(\lambda,\mu,z,\bar z)f^{(+)}(\mu,z,\bar z)=
O\left(\frac{1}{\lambda}\right)e^{\lambda z},
\label{residue_+-}
\eeq
as $\lambda\rightarrow\infty_-$.

Similarly if $f^{(-)}(\lambda,z,\bar z)$ is the 
following formal series in $\lambda$
\beq
f^{(-)}(\lambda,z,\bar z)=\left\{
  \sum_{k=-N}^{\infty}\frac{f^{(-)}_k(z,\bar
    z)}{\lambda^k}\right\}e^{\lambda z},
\label{formal_-}
\eeq
then
\beq 
\left.2\pi i\,\hbox{\rm res}\vphantom{\int}\right|_{\mu=\infty_-}
\!\!\!\!\!\!\!\!\omega(\lambda,\mu,z,\bar z)f^{(-)}(\mu,z,\bar z)= 
\left\{
  \sum_{k=-N}^{0}\frac{f^{(-)}_k(z,\bar
    z)}{\lambda^k}+O\left(\frac{1}{\lambda}\right)\right\}e^{\lambda z}, 
\label{residue_--}
\eeq
as $\lambda\rightarrow\infty_-$ and 
\beq
\left.2\pi i\,\hbox{\rm res}\vphantom{\int}\right|_{\mu=\infty_-}
\!\!\!\!\!\!\!\!\omega(\lambda,\mu,z,\bar z)f^{(-)}(\mu,z,\bar z)=
O\left(\frac{1}{\lambda}\right)e^{\lambda\bar z},
\label{residue_-+}
\eeq
as $\lambda\rightarrow\infty_+$.
\end{enumerate}
\end{statement}

\begin{remark}
The formulas (\ref{residue_++}), (\ref{residue_+-}), (\ref{residue_--}),
(\ref{residue_-+}) require some comments. On the left-hand side we
have formal series in $\mu$, thus the analytic definition of the residue as
an integral does not work. Fortunately the terms, containing
exponents in $\mu$ annihilate each other and we can define the
residue as the coefficient of $1/\mu$:
\beq
\left.\hbox{\rm res}\vphantom{\int}\right|_{\mu=\infty}\sum_{k=-N}^{\infty}
\frac{a_k}{\mu^k} d\mu=a_1
\eeq
\end{remark}
Proof of Lemma \ref{lemma:cba}. 
\begin{enumerate}
\item We have constructed $\tilde\omega(\lambda,\mu,z,\bar z)$ as a
meromorphic function. The function $<\tilde\omega_x(\mu,\mu,z,\bar z)>_x$ 
is meromorphic in $\mu$ on the whole $\Gamma$ and $dp_1(\mu)$ is a
meromorphic 1-form on $\Gamma$. Thus (\ref{cba}) gives us a
meromorphic function with appropriate tensor properties. 
\item From Lemma \ref{lemma:asymptotic} it 
follows that for generic $\mu$ 
\begin{enumerate}
\item $\vec\psi(\nu,z,\bar z)$ is nonsingular at the points 
$\nu=\mu$, $\nu=\sigma\mu$, 
\item $<\omega_x(\mu,\mu,z,\bar z)>_x\ne0$ 
\item $w_1(\lambda)w_1^{-1}(\mu)=w_2(\lambda)w_2^{-1}(\mu)=1$ iff 
$\lambda=\mu$.
\end{enumerate}
If $\mu$ fulfill these conditions, if $\lambda\ne\mu$ and if
$\vec\psi(\lambda,z,\bar z)$ is nonsingular, then formulas
(\ref{cba_reg_2}), (\ref{cba_reg_3}), (\ref{cba}), define a
nonsingular function of $z$, $\bar z$. 
\item Let $\lambda\rightarrow\mu$. To calculate the asymptotics of the
CBA kernel we expand the function $w_1(\lambda)$ in (\ref{cba_reg})
in a series in $\lambda-\mu$ using (\ref{p's}), (\ref{dp's}). We see,
that we have a first-order pole and the normalization in (\ref{cba})
is chosen so that the residue is exactly 1.
\item To calculate the asymptotics of $\tilde\omega(\lambda,\mu,z,\bar z)$
as $\mu\rightarrow\infty$ we substitute the asymptotic expansion
(\ref{phias}) to (\ref{cba_1}) and use the following asymptotic formula:
\beq
\int^{z} e^{\lambda\bar z}f(\lambda,z',\bar z') dz'+
e^{\lambda\bar z}g(\lambda,z',\bar z') d\bar z'=
e^{\lambda\bar z}\sum_{k=1}^{\infty}(-1)^{k-1}\frac{\partial_{\bar z}^{k-1}
g(\lambda,z,\bar z)}{\lambda^{k}}, 
\label{as_integration}
\eeq
where the functions $f(\lambda,z',\bar z')$, $g(\lambda,z',\bar z')$
are some asymptotic series in $\lambda$
\beq
f(\lambda,z',\bar z')=\sum_{k=0}^{\infty}\frac{f_k(z,\bar z)}{\lambda^k}, \ 
g(\lambda,z',\bar z')=\sum_{k=0}^{\infty}\frac{g_k(z,\bar z)}{\lambda^k} 
\label{as_integrands}
\eeq
such that
\beq
\partial_{\bar z}\left( e^{\lambda\bar z}f(\lambda,z',\bar z')\right)=
\partial_{z}\left( e^{\lambda\bar z}g(\lambda,z',\bar z')\right),
\eeq
and all expansion coefficients $f_k(z,\bar{z})$, $g_k(z,\bar{z})$
are smooth functions, all derivatives of these functions are bounded
in the $z$-plane.

The proof of the formula (\ref{as_integration}) is rather standard and
we do not reproduce it here. 
\item The last statement of the Lemma follows directly from the
transformation rules for the Bloch functions.
\end{enumerate}
Proof of Statement~\ref{statement:cba}.
Assume for definiteness, that $\lambda\rightarrow\infty_+$. 
If $\mu\rightarrow\infty_+$ then:
\beq
\omega(\lambda,\mu,z,\bar{z})=\left[ \frac{1}{2\pi i(\mu-\lambda)}+
\sum\limits_{i,j>0}
\frac{\omega^{(++)}_{ij}(z,\bar{z})}{\lambda^i\mu^j}
\right]
e^{(\lambda-\mu)\bar{z}}d\mu.
\label{omega_at_infty}
\eeq
From (\ref{omega_at_infty}) it follows, that 
$$
\left.2\pi i\,\hbox{\rm res}\vphantom{\int}\right|_{\mu=\infty_+}
\!\!\!\!\!\!\!\!\omega(\lambda,\mu,z,\bar z)f^{(+)}(\mu,z,\bar z)=
$$
\beq 
=\left.\hbox{\rm res}\vphantom{\int}\right|_{\mu=\infty_+}
\frac{d\mu}{\mu-\lambda}e^{(\lambda-\mu)\bar{z}}f^{(+)}(\mu,z,\bar z)+
O\left(\frac{1}{\lambda}\right)e^{\lambda\bar z}. 
\label{residue_calc_1}
\eeq
Applying the well-known formula:
\beq
\left.\hbox{\rm res}\vphantom{\int}\right|_{\mu=\infty}
\frac{d\mu}{\mu-\lambda}\mu^n=\left\{\barr{ll} \lambda^n & n\ge0 \\ 0 & n<0
\earr\right.
\label{Cauchy}
\eeq
we get the right-hand side of (\ref{residue_++}).

If $\mu\rightarrow\infty_-$, then using (\ref{as_omega_+}) we get:
\beq
\omega(\lambda,\mu,z,\bar{z})=\left[\sum\limits_{i,j>0}
\frac{\omega^{(+-)}_{ij}(z,\bar{z})}{\lambda^i\mu^j}
\right]e^{\lambda\bar{z}-\mu z}d\mu.
\label{different_inft's}
\eeq
Formula (\ref{residue_-+}) follows automatically from 
from (\ref{different_inft's}). 

Formulas (\ref{residue_--}), (\ref{residue_+-}) are proved exactly
in the same way.
\begin{remark}
A formula, representing the Cauchy kernel on a Riemann surface as a
semi-infinite sum of quadratic combinations of eigenfunctions first
arose in the article \cite{KN87} by I.~M.~Krichever and S.~P.~Novikov. In 
\cite{KN87} the spatial variable $x$ was discrete. In \cite{GO89}
the spatial variable was continuous, and the CBA kernel was defined as
an integral similar to (\ref{cba_1}).
\end{remark}
\begin{remark}
A form, similar to $\tilde\omega(\lambda,\mu,z,\bar z)$, but with 
integration path, starting from a finite point $Z$ in the $z$-plane, 
arises in the theory of finite Darboux transformations (see
\cite{MS91}, formula (6.1.18)). In fact, the same object arose in the
Generalized Weierstrass map. We do not discuss this analogy in our
text, but it seems possible that this analogy has some deep roots.
\end{remark}

\subsection{Modified Novikov-Veselov equations.}
\label{sec:mnv}

\begin{lemma}
\label{lemma:deform_1}
Let $\mu$ be a generic point in $\Gamma$. Define a deformation
$\delta^{(\mu)}$ of the function $\vec\psi(\lambda,z,\bar z)$ 
associated with this point by:
\beq
\barr{c}
\delta^{(\mu)} \psi_1(\lambda,z,\bar z)=
\tilde\omega(\lambda,\mu,z,\bar z)\psi_1(\mu,z,\bar z)+
\alpha(\lambda,\mu) \psi_1(\lambda,z,\bar z) 
\\
\delta^{(\mu)} \psi_2(\lambda,z,\bar z)=
\tilde\omega(\lambda,\mu,z,\bar z)\psi_2(\mu,z,\bar z)+
\alpha(\lambda,\mu) \psi_2(\lambda,z,\bar z),
\earr
\label{delta_psi_1}
\eeq
where $\alpha(\lambda,\mu)$ is a meromorphic function in $\lambda$,
uniquely fixed by the following requirement
\beq
\left.\delta^{(\mu)}\psi_1(\lambda,z,\bar z)+
\delta^{(\mu)}\psi_2(\lambda,z,\bar z)\right|_{z=z_1}=0.
\label{norm_defs}
\eeq
Then (\ref{delta_psi_1}) is a deformation of the Bloch function,
corresponding to the following deformation of the Dirac operator
\beq
\delta^{(\mu)}L=\left[\barr{cc}0 & -\delta^{(\mu)}_1 U(z,\bar z) \\ 
\delta^{(\mu)}_2 U(z,\bar z) & 0 \earr\right],
\label{delta_L_1}
\eeq
where
\beq
\barr{c}
\delta^{(\mu)}_1 U(z,\bar z)=\psi_1(\mu,z,\bar z)\psi_2(\sigma\mu,z,\bar z) 
\\
\delta^{(\mu)}_2 U(z,\bar z)=\psi_2(\mu,z,\bar z)\psi_1(\sigma\mu,z,\bar z).
\earr
\label{delta_U_1}
\eeq

For generic $\mu$ these deformations result in non-self-adjoint Dirac 
operators ($\delta^{(\mu)}_1 U(z,\bar z)\ne\delta^{(\mu)}_2 U(z,\bar z)$)
with complex-valued potentials. But the
the Bloch function and the Bloch variety are well defined for
such Dirac operators.
\end{lemma}
Proof of Lemma \ref{lemma:deform_1}. A simple direct calculation shows,
that 
\beq
\left(L+\epsilon\delta^{(\mu)}L\right) 
\left(\vec\psi(\lambda,z,\bar z)+
\epsilon\delta^{(\mu)}\vec\psi(\lambda,z,\bar z)\right)=O(\epsilon^2),
\label{deformed_dirac}
\eeq
Thus (\ref{delta_psi_1}) defines deformations of the eigenfunctions,
corresponding to (\ref{delta_L_1}). From (\ref{periods_omega}) it
follows that the new eigenfunctions satisfy (\ref{blochsol}) with the
same multipliers $w_1(\lambda)$, $w_2(\lambda)$ as the old ones. Thus
the new eigenfunctions are defined on the same curve $\Gamma$ and have
the same periodicity properties. The last property plays the key
role, when we prove below that 
the functionals $h_{2k+1}[U]$ are invariant under some deformations.

The deformations of the Dirac operator, generated by all
$\delta^{(\mu)}$ form a linear space. Deformations preserving the
class of self-adjoint Dirac operators with real potentials are the
most interesting ones. Let us check that the subspace of such
deformations is sufficiently large.

\begin{lemma}
\label{lemma:deform_4}
Denote by $\Delta^{(\mu)}$ the following linear combination of
deformations $\delta^{(\mu)}$:
\beq
\Delta^{(\mu)}=
\frac{\delta^{(\mu)}+\delta^{(\sigma\mu)}}
{<\tilde\omega_x(\mu,\mu,z,\bar z)>_x}+
\frac{\delta^{(\theta\mu)}+\delta^{(\sigma\theta\mu)}}
{<\tilde\omega_x(\theta\mu,\theta\mu,z,\bar z)>_x}
\label{delta_4}
\eeq
(recall that $<\tilde\omega_x(\mu,\mu,z,\bar z)>_x$ is an even
function in $\mu$ (see (\ref{even})) and does not depend on $z$, $\bar z$.) 

Then $\Delta^{(\mu)}$ acts on the space of self-adjoint Dirac
operators with real potentials, i.e.
\beq
\barr{c}
\Delta^{(\mu)}_1 U(z,\bar z)=\Delta^{(\mu)}_2 U(z,\bar z)=
\overline{\Delta^{(\mu)}_1 U(z,\bar z)}=
\overline{\Delta^{(\mu)}_2 U(z,\bar z)}=
\\ \\
=\frac{\psi_1(\mu,z,\bar z)\psi_2(\sigma\mu,z,\bar z)+
\psi_1(\sigma\mu,z,\bar z)\psi_2(\mu,z,\bar z)}
{<\tilde\omega_x(\mu,\mu,z,\bar z)>_x}+
\frac{\psi_1(\theta\mu,z,\bar z)\psi_2(\sigma\theta\mu,z,\bar z)+
 \psi_1(\sigma\theta\mu,z,\bar z)\psi_2(\theta\mu,z,\bar z)}
{<\tilde\omega_x(\theta\mu,\theta\mu,z,\bar z)>_x}
\label{delta_U_4}
\earr
\eeq
and the Bloch variety $\Gamma[U]$ defined in Section~\ref{sec:bloch} is 
invariant under these deformations. 
\end{lemma}
Proof of Lemma \ref{lemma:deform_4}. Formula (\ref{delta_U_4}) follows
directly from (\ref{delta_U_1}) and (\ref{delta_4}). The Bloch
function $\vec\psi(\lambda,z,\bar z)$ for real $U(z,\bar z)$ has
symmetry property (\ref{n}), thus the right-hand side of
(\ref{delta_U_4}) is real. As we pointed out in the previous Lemma,
all deformations generated by $\delta^{(\mu)}$ do not change the
Bloch variety $\Gamma[U]$ and  the multipliers $w_1$, $w_2$. This
completes the proof. If Statement~\ref{st2} is fulfilled, then these
flows do not change $h_{2k+1}[U]$.

Equations (\ref{delta_U_4}) are essentially nonlocal and rather
complicated. Namely the right-hand side is expressed in terms of the
Bloch function, and it is difficult to calculate Bloch solutions
either analytically or numerically. Fortunately the space of
deformations generated by all $\Delta^{(\mu)}$ contains simpler
equations such that the right-hand side can be expressed via
$U(z,\bar z)$ in terms of quadratures.

Let $\mu\rightarrow\infty$. If Statement~\ref{st2} is fulfilled we may
expand $\Delta^{(\mu)}$ to the following asymptotic series
\beq
\Delta^{(\mu)}U(z,\bar z)=
-2\frac{d\mu}{idp(\mu)}\sum_{k=0}^{\infty}\frac{K_{2k+1}[U]}{\mu^{2k+2}}
+\frac{\overline{K_{2k+1}[U]}}{\bar\mu^{2k+2}}.
\label{novikov-veselov_1}
\eeq
(We write the term $d\mu/dp(\mu)$ to gain some standard normalization
of the MNV hierarchy. If we omit this multiplier, our
expansion coefficients will be linear combinations of Novikov-Veselov
generators with constant coefficients, which in most situations is not
essential). 

Any $K_{2k+1}[U]$ is a quadratic polynomial of $\phi_l(z,\bar z)$,
$\chi_l(z,\bar z)$, $\l=1,\ldots,2k$, with coefficients depending on the 
$h_{2l-1}$ and $c_{2l}$, $l=1,\ldots,k$, where the $c_{2l}$ are
coefficients of the asymptotic expansion (\ref{even_2}).

For any odd integer $2k+1$, $k\ge0$ we have the
following pair of flows on the space of real double-periodic functions:
\beq
\frac{\partial U(z,\bar z, t_{2k+1})}{\partial t_{2k+1}}=
2\,\hbox{Re}\,K_{2k+1}[U].
\label{novikov_veselov_2}
\eeq
\beq
\frac{\partial U(z,\bar z, \tilde t_{2k+1})}{\partial \tilde t_{2k+1}}=
2\,\hbox{Im}\,K_{2k+1}[U].
\label{novikov_veselov_3}
\eeq

\begin{definition}
The equations (\ref{novikov_veselov_2}), (\ref{novikov_veselov_3}) are
called the {\bf Modified Novikov-Ve\-se\-lov equations} (MNV).
\end{definition}

\begin{remark}
To define the MNV flows associated to (\ref{novikov_veselov_2}),
(\ref{novikov_veselov_3}) 
it is sufficient to have a formal expansion for the Bloch function and
the function $h(\lambda)$. Thus these flows are well-defined for any
smooth potential. But without Statement~\ref{st2} we can not prove
that they conserve the functional $h_{2k+1}[U]$.
\end{remark}

Using this definition of the MNV hierarchy we may immediately prove 
Statement~\ref{lemma:conservation}. All functionals $h_{2k+1}[U]$ defined
in Section \ref{sec:bloch} are integrals of motion for all flows 
$\Delta^{(\mu)}$, thus they are conservation laws for their expansion
coefficients. 

The representation (\ref{novikov_veselov_2}) may look rather unusual. To
check that our definition of the MNV equations coincides with the
standard one, let us calculate the deformation of the Bloch function,
corresponding to the flows (\ref{novikov_veselov_2}). To gain a standard
answer we shall use a normalization of the Bloch function different
from th eone used 
earlier. Instead of (\ref{norm_defs}) we assume that the function
$\alpha(\lambda,\mu)$ in (\ref{delta_psi_1}) is identically zero. 

\begin{statement}
\label{lemma:standard_mnv}
Let the Statement~\ref{st2} be valid. Then 
the deformation of the Bloch function corresponding to the flows 
(\ref{novikov_veselov_2}) has the following representations:
$$
\frac{\partial\vec\psi(\lambda,z,\bar z,t_{2k+1})}
{\partial t_{2k+1}}=
$$
\beq
=2\pi i\left\{\left.\hbox{\rm res}\vphantom{\int}\right|_{\mu=\infty_+}
\!\!\!\!\!\!\!\!\!\!\!
+\!\left.\hbox{\rm res}\vphantom{\int}\right|_{\mu=\infty_-}\right\}
\omega(\lambda,\mu,z,\bar{z},t_{2k+1},)
\mu^{2k+1}\vec\psi(\mu,z,\bar z,t_{2k+1}),
\label{novikov_veselov_2_psi_1}
\eeq
\beq
=\lambda^{2k+1}\vec\psi(\lambda,z,\bar z,t_{2k+1})+
\left\{\barr{l} 
O\left(\frac{1}{\lambda}\right)e^{\lambda \bar z}\ \hbox{as} 
\ \lambda\rightarrow\infty_+
\\
O\left(\frac{1}{\lambda}\right)e^{\lambda z} \ \hbox{as}
\ \lambda\rightarrow\infty_-
\earr\right.
\label{novikov_veselov_2_psi_2}
\eeq
\beq
=\left\{\partial_z^{2k+1}+\partial_{\bar z}^{2k+1}+
\sum_{l=0}^{2k}W_l(z,\bar z)\partial_z^{l}+
\overline{W_l(z,\bar z)}\partial_{\bar z}^{l}\right\}
\vec\psi(\lambda,z,\bar z,t_{2k+1}),
\label{novikov_veselov_2_psi_3}
\eeq
\end{statement}
Proof of Statement~\ref{lemma:standard_mnv}, given the validity of
Statement~\ref{st2}. By definition 
$$
-\frac{idp(\mu)}{2}\Delta^{(\mu)} \vec\psi(\lambda,z,\bar z)=
$$
\beq
\barr{c}
=\pi i\left[\omega(\lambda,\mu,z,\bar z)\vec\psi(\mu,z,\bar z)+
[\omega(\lambda,-\mu,z,\bar z)\vec\psi(-\mu,z,\bar z)\right]_
{\mu\rightarrow\infty_+}+
\\
+\pi i \left[\omega(\lambda,\bar\mu,z,\bar z)\vec\psi(\bar\mu,z,\bar z)+
[\omega(\lambda,-\bar\mu,z,\bar z)\vec\psi(-\bar\mu,z,\bar z)\right]_
{\mu\rightarrow\infty_-}
\earr
\label{Delta_psi}
\eeq
(in this formula we use $\sigma\mu=-\mu$, $\theta\mu=-\bar\mu$).

The residues in (\ref{novikov_veselov_2_psi_1}) are exactly the
coefficients of the terms $d\mu/\mu^{2k+2}$ and $d\bar
\mu/\bar\mu^{2k+2}$ respectively. Thus these coefficients gives us the
action of $K_{2k+1}[U]$ and $\overline{K_{2k+1}[U]}$ on the Bloch
function.

Formula (\ref{novikov_veselov_2_psi_2}) follows directly from 
(\ref{novikov_veselov_2_psi_1}) and formulas
(\ref{residue_++})-(\ref{residue_-+}).

To prove (\ref{novikov_veselov_2_psi_3}) let us substitute the asymptotic
expansions (\ref{as_omega_+}), (\ref{as_omega_-}) in
(\ref{novikov_veselov_2_psi_1}). We get:
$$
\frac{\partial\psi_1(\lambda,z,\bar z,t_{2k+1})}
{\partial t_{2k+1}}=
$$
\beq
\barr{l}
=\left\{R^{(+)}_{2k+1}[\partial_{\bar z}]+\sum\limits_{l=0}^{2k}
\phi_{2k+1-l}(z,\bar z) R^{(+)}_{l}[\partial_{\bar z}]\right\}
\psi_1(\lambda,z,\bar z,t_{2k+1})+
\\
+\left\{\sum\limits_{l=0}^{2k}
\chi^{(-)}_{2k+1-l}(z,\bar z) R^{(-)}_{l}[\partial_{z}]\right\}
\psi_2(\lambda,z,\bar z,t_{2k+1})=
\earr
\eeq
\beq
=\partial_{\bar z}^{2k+1}\psi_1(\lambda,z,\bar z,t_{2k+1})+
U(z,\bar z)\partial_{z}^{2k}\psi_2(\lambda,z,\bar z,t_{2k+1})+
\,\hbox{lower order terms}=
\eeq
\beq
=\left\{\partial_{\bar z}^{2k+1}+ \partial_{z}^{2k+1} \right\}
\psi_1(\lambda,z,\bar z,t_{2k+1})+
\,\hbox{lower order terms},
\eeq
(the $\chi^{(-)}_{l}(z,\bar z)$ denote the  expansion coefficients of the
function $\psi_1(\lambda,z,\bar z)$ at the point $\lambda=\infty_-$). 
A similar calculation shows the following relation:
\beq
\frac{\partial\psi_2(\lambda,z,\bar z,t_{2k+1})}
{\partial t_{2k+1}}=
\left\{\partial_{\bar z}^{2k+1}+ \partial_{z}^{2k+1} \right\}
\psi_2(\lambda,z,\bar z,t_{2k+1})+
\,\hbox{lower order terms}. 
\eeq
Statement \ref{lemma:standard_mnv} is proved.

From (\ref{novikov_veselov_2_psi_3}) it follows that the MNV equations
are the compatibility conditions for a pair of differential operators on
the space of zero eigenfunctions of $L$. It is well-known in soliton
theory that such compatibility conditions are equivalent to the 
existence of standard $L-A-B$ representations (see \cite{Bog87}).

\begin{remark}
In this article we constructed a soliton hierarchy in terms of the
Cauchy-Baker-Akhiezer kernel.

Using the CBA kernel we construct, in fact, a 
much wider hierarchy, including essentially nonlocal equations. All these
equations preserve the spectral curve. In Section \ref{sec:conf} we
show that in fact the deformations corresponding to conformal
transformations of the Euclidean space, lie in this wider hierarchy. 

In terms of $U(z,\bar{z})$ this hierarchy looks rather unnatural. But
we may treat it as a system of differential equations on a bigger
collection of functions simultaneously; namely we may consider the potential
$U(z,\bar{z})$ and the wave function in a finite number of fixed points
on the spectral curve as unknown functions, connected by the Dirac equations. 
Similar systems associated with one-dimensional soliton equations were
discussed in the literature (see \cite{Mel90} and references therein)
from both a mathematical and a physical point of view. In \cite{Or91} it
was shown, that starting from the KP equation we get a hierarchy 
naturally containing many other well-known soliton systems.
\end{remark}

\end{document}